# Considering Multiple Uncertainties in Stochastic Security-Constrained Unit Commitment Using Point Estimation Method

Mahdi Mehrtash*, Mahdi Raoofat, Mohammad Mohammadi, Mohammad Hossein Zakernejad
Department of Power and Control Eng.
Shiraz University
Shiraz, Iran
*eng.mehrtash@yahoo.com

Hamidreza Zareipour
Department of Electrical and Computer Eng.
University of Calgary
Calgary, Canada
h.zareipour@ucalgary.ca

*Abstract*— Security-Constrained Unit Commitment (SCUC) is one of the most significant problems in secure and optimal operation of modern electricity markets. New sources of uncertainties such as wind speed volatility and price-sensitive loads impose additional challenges to this large-scale problem. This paper proposes a new Stochastic SCUC using point estimation method to model the power system uncertainties more efficiently. Conventional scenario-based Stochastic SCUC approaches consider the Mont Carlo method; which presents additional computational burdens to this large-scale problem. In this paper we use point estimation instead of scenario generating to detract computational burdens of the problem. The proposed approach is implemented on a six-bus system and on a modified IEEE 118-bus system with 94 uncertain variables. The efficacy of proposed algorithm is confirmed, especially in the last case with notable reduction in computational burden without considerable loss of precision.

*Keywords— Bender's Decomposition; Two-Point Estimation method; scenario-based approach; Stochastic SCUC*

## I. Introduction

Nowadays, penetration level of renewable energy resources is increasing rapidly. It is reported by the Global Wind Energy Council (GWEC) that global wind energy installations increased by 35,282 MW in 2013, about 13 percent of all installations until 2012 [1].

Renewable energy resources like wind are highly probabilistic and intensify uncertainties of planning and operation of electric power systems. Modern wind forecasting techniques can predict total power generation of wind farms for next few hours with an error of 5% to 20% [2], [3]. Load variation is another source of uncertainty that is worsening with emersion of Demand Response (DR) and price-sensitive loads. In this uncertain environment, secure-optimal operation of markets is a vital responsibility of ISO and market operators, therefore Security Constrained Unit Commitment (SCUC) [4] and Stochastic SCUC [5] have been studied and reported widely in literature.

Due to increase of uncertainties of power systems, Stochastic SCUC has become more necessary, recently. Random behavior of Plug-in Electric Vehicles (PEVs) besides the stochastic wind power generation is considered in stochastic SCUC in [6]. The authors of [7] propose a Stochastic SCUC model for hourly coordination of wind power and pumped-storage hydro generations in order to increase the dispatchability of wind energy.

Two distinct approaches have been proposed in literature for Stochastic SCUC: scenario-based and interval optimization approaches [8].

In scenario-based approaches, sufficient number of probable scenarios are generated and the problem is solved for each scenario. Then, output variables are deduced from the set of individual results. Monte-Carlo Simulation (MCS) is one of the most famous techniques in this category and is very flexible in modeling any kind of uncertainty in any complex problem [9]. In another technique, first the Probability Distribution Functions (PDFs) of random variables are quantized and then, all combinations of quantized values are studied scenario by scenario. Reference [10] presents a scenario-based SCUC model considering wind power intermittency and volatility.

In interval optimization approach, the probable interval of objective function is obtained with a significantly less computation burden. The method is utilized for stochastic SCUC in [8], and the results are compared with scenario-based approach. The authors confirm that scenario-based approach is very time-consuming with a high computational burden, but with a more stability in results. In addition, interval optimization may not be applicable for problems with discrete uncertain variables. Necessity to work on speeding up the scenario-based approach as a better choice is emphasized in [8].

Point Estimation (PE) is a well-known analytical method in which each PDF is concentrated into few points, provided by its first central moments. All probable combinations of these concentrations constitute the set of scenarios, then the calculation speed is improved significantly, especially in large-scale problems. PE method has been used in literature for different power system problems such as probabilistic power flow [11], probabilistic optimal load flow [12], small signal stability analysis [13], etc. PE method has not been reported for stochastic SCUC, so far.

In this paper, scenario-based Stochastic SCUC model with uncertain loads and wind generations is formulated and solved using PE method, which results in speeding up the calculation without significant loss of accuracy. Numerical studies on 6-bus system, which includes one uncertain wind generator and three random loads is illustrated in detail as the first case study. In the second case, the algorithm is implemented on a modified 118-bus system with three wind generators and 91 stochastic loads with noticeable reduction in calculation time without considerable decrease of accuracy. Solving such a large scenario-based stochastic SCUC example is not reported in literature, so far.

The rest of the paper is organized as follows. Next section presents PE method in stochastic analysis. Section III and IV explain the Stochastic SCUC formulation and the proposed Stochastic SCUC with PE method, respectively. Numerical examples are presented in section V of the paper. At last, a conclusion is presented in section VI.

## II. POINT ESTIMATION METHODS

Point Estimation (PE) methods concentrate the statistical information of each random variable into some points, named *concentrations*, provided by the first few central moments of the variables [14]. For each random input variable, the $k^{th}$ concentration ($P_{l,k}$, $W_{l,k}$) can be defined as a pair of *location* $P_{l,k}$ and *weight* $W_{l,k}$. Equation (1) determines the location $P_{l,k}$. Where $\mu_{pl}$ and $\sigma_{pl}$ are the mean and standard deviation of input random variable $Pl$, and $\xi_{l,k}$ is its standard location. For TPE method the standard location $\xi_{l,k}$ and the weight $W_{l,k}$ can be obtained using equations (2) and (3). Where $m$ is the number of input random variables, and $\lambda_{l,j}$ denotes the *jth* standard central moment of the random variable $Pl$. After determining the concentrations, there are different approaches to solve the main probabilistic problem. For problems with normal PDFs or with insignificant asymmetry in their PDFs Two-Point Estimation (TPE) method is accurate and common in literature. One of the most accepted approaches for TPE proposes to analyze the problem in $2\times m$ scenarios [15]. Once all the concentrations ($P_{l,k}$, $W_{l,k}$) are calculated, each scenario is a deterministic problem evaluated in the point ($\mu_{p1},\mu_{p2},…,P_{l,k},…,\mu_{pm}$) consists of one concentration of variable $Pl$ and mean values of other variables. $Z$ is the vector of output random variables, and $Z_{(l,k)}$ is the output of corresponding scenario. Lastly, (4) determines the *jth* raw moment of the output random variables. Where, $m$ is the number of random variables, and $C$ is the number of concentrations assumed for each random variable.

$$p_{l,k} = \mu_{pl} + \xi_{l,k}\, \sigma_{pl} \quad (1)$$

$$\xi_{l,1} = \frac{\lambda_{l,3}}{2} + \sqrt{m + \left(\frac{\lambda_{l,3}}{2}\right)^2} \;,\; \xi_{l,2} = \frac{\lambda_{l,3}}{2} - \sqrt{m + \left(\frac{\lambda_{l,3}}{2}\right)^2} \quad (2)$$

$$w_{l,1} = -\frac{1}{m}\frac{\xi_{l,2}}{\xi_{l,1}-\xi_{l,2}} \;,\; w_{l,2} = \frac{1}{m}\frac{\xi_{l,1}}{\xi_{l,1}-\xi_{l,2}} \quad (3)$$

$$\mu_j = E[Z^j] \cong \sum_{l=1}^{m}\sum_{k=1}^{C} w_{l,k}\, (Z(l,k))^j \quad (4)$$

## III. STOCHASTIC SCUC DEFINITION AND FORMULATION

### A. Stochastic SCUC

In deterministic SCUC, the problem is solved just with certain variables, which is shown not to be efficient in presence of stochastic variables like wind generation or loads [16]. To overcome this inefficiency, some different algorithms are proposed in literature, which can be categorized in scenario-based and interval optimization algorithms. The former category presents better information about the system and is more stable in convergence [8].

Here, a new scenario-based approach is proposed for SCUC which uses TPE method to overcome the computational budrn of scenario-based approaches without significant loss of accuracy. As Fig. 1 shows, an initial optimal operating point, X0, is obtained using deterministic SCUC with the Most Probable Values (MPVs) for random variables. This case is called the base case, here. As Fig.1 shows, in different possible real-time cases, the feasible region of operation which are shown by S1, S2, ... may be unreachable from X0, due to system or unit restrictions. The final optimum operating point for the next day, X*, should be obtained so that in any probable case in real time, the system can transfer to a new optimal point with a corrective action. Transferring from X0 to X* in the time of scheduling is called preventive action, and moving in real time from X* to actual optimal operating point $X_i$ of $S_i$ is corrective action [10]. In proposed algorithm, preventive actions are administered using Bender's cuts.

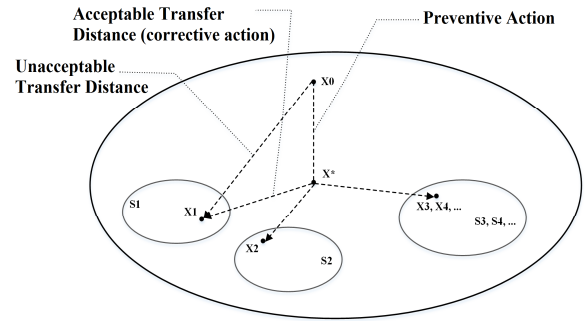

Fig. 1. Corrective and preventive actions in Stochastic SCUC problem

## B. Problem Formulation

The objective function (5) consists of generation cost; and startup and shutdown costs of thermal units.

$$\min \sum_{i=1}^{NG} \sum_{t=1}^{NT} [F_{ci}(P_{it}) \cdot I_{it} + SU_{it} + SD_{it}] \quad (5)$$

Base case is run as a deterministic SCUC with MPVs of each random variable. The base case constraints include system power balance as presented by (6), system spinning and operating reserves as (7), ramping up/down limits as (8), minimum ON/OFF time limits as (9), unit generation limits as (10) and transmission lines constraints as (11).

$$\sum_{i=1}^{NG} P_{it} \cdot I_{it} + \sum_{i=1}^{NW} P_{w,it}^f = P_{D,t}^f \quad (t=1,\dots,NT) \quad (6)$$

$$\sum_{i=1}^{NG} R_{S,it} \cdot I_{it} \geq R_{S,t} \quad (t=1,\dots,NT)$$

$$\sum_{i=1}^{NG} R_{O,it} \cdot I_{it} \geq R_{O,t} \quad (t=1,\dots,NT) \quad (7)$$

$$P_{it} - P_{i(t-1)} \leq [1 - I_{it}(1 - I_{i(t-1)})] \cdot UR_i + I_{it}(1 - I_{i(t-1)}) \cdot P_i^{min}$$
$$(t=1,\dots,NT; i=1,\dots,NG)$$
$$P_{i(t-1)} - P_{it} \leq [1 - I_{i(t-1)}(1 - I_{it})] \cdot DR_i + I_{i(t-1)}(1 - I_{it}) \cdot P_i^{min}$$
$$(t=1,\dots,NT; i=1,\dots,NG) \quad (8)$$

$$[X_{i(t-1)}^{on} - T_i^{on}] \cdot [I_{i(t-1)} - I_{it}] \geq 0 \quad (t=1,\dots,NT; i=1,\dots,NG)$$
$$[X_{i(t-1)}^{off} - T_i^{off}] \cdot [I_{it} - I_{i(t-1)}] \geq 0 \quad (t=1,\dots,NT; i=1,\dots,NG) \quad (9)$$

$$P_i^{min} \cdot I_{it} \leq P_{it} \leq P_i^{max} \cdot I_{it} \quad (t=1,\dots,NT; i=1,\dots,NG) \quad (10)$$

$$-PL^{max} \leq SF \cdot (K_P \cdot P_t - K_D \cdot P_{Dt}) \leq PL^{max} \quad (11)$$

Constraints for each scenario are approximately similar with the base case with small differences. The constraints are system power balance equation (12), unit generation limits (13), and dispatch adjustment capabilities of generating units (14). Transmission lines constraints (15) are also met in each scenario to guarantee the network security.

$$\sum_{i=1}^{NG} P_{it}^s \cdot I_{it} + \sum_{i=1}^{NW} P_{w,it}^s = \sum_{b=1}^{NB} P_{D,bt}^s \quad (12)$$

$$P_i^{min} \cdot I_{it} \leq P_{it}^s \leq P_i^{max} \cdot I_{it} \quad (13)$$

$$P_{it}^s - P_{it} \leq R_i^{up} \cdot I_{it}$$
$$P_{it} - P_{it}^s \leq R_i^{dn} \cdot I_{it} \quad (14)$$

$$-PL^{max} \leq SF \cdot (K_P \cdot P_t^s - K_D \cdot P_{Dt}^s) \leq PL^{max} \quad (15)$$

## IV. SOLUTION METHODOLOGY

Scenario-based Stochastic SCUC problem modeled in (5)–(15) is a non-convex, large scale, and non-deterministic polynomial-time hard (NP-hard) problem. Solving this problem for large-scale power systems would be intractable without decomposition [17]. The Bender's Decomposition (BD) technique is implemented to decompose the problem into a master problem and several sub-problems.

Fig. 2 shows the flowchart of the proposed method for Stochastic SCUC problem. According to flowchart of Fig. 2, in proposed method the Stochastic SCUC problem is solved in three steps using BD. Firstly, the master Unit Commitment (UC) problem is administered considering all deterministic constraints of the units. Then, the network constraints are evaluated and met to fulfill deterministic SCUC. The results of these two steps should be evaluated in scenario check sub-problem, for each scenario. Every mismatch or violation will be referred to the first steps by BD technique, so that finally all constraints are met in all probable scenarios. This paper proposes that instead of random scenarios proposed in [8] or scenarios provided by LHS in [10], TPE is used to generate efficient scenarios to be checked in the third step of the algorithm. The efficacy, speed and accuracy of this method is investigated by numerical studies.

### A. Master Unit Commitment Problem

In order to minimize the operating cost of the base case, the objective function (5) is optimized with constraints of (6)–(10) as a deterministic mixed-integer programming. In this problem all uncertain input variables are fixed to their MPVs.

### B. Hourly Network Check for the Base Case

After obtaining the results of the master UC problem, the transmission network constraints of (11) are evaluated and met using BD technique by (16) and (17) [8]. $S_t$ is the objective value, and if its minimum value remains larger than a predefined threshold, it means at least one network security constraint is violated. In this situation a feasibility cut (17) is generated and used in the next iteration of the master UC problem. If there is no any violation, it shows that a feasible and optimal solution has been found for the deterministic SCUC problem.

$$\min s_t$$
$$s.t.$$
$$-1 \cdot s_t \leq PL^{max} - SF \cdot (K_P \cdot \hat{P}_t - K_D \cdot P_{Dt}) \quad ; \quad \lambda_{1,t}$$
$$-1 \cdot s_t \leq PL^{max} + SF \cdot (K_P \cdot \hat{P}_t - K_D \cdot P_{Dt}) \quad ; \quad \lambda_{2,t}$$
$$s_t \geq 0 \quad (16)$$

$$-(\lambda_{1,t} - \lambda_{2,t})^T \cdot K_P \cdot (P_t - \hat{P}_t) + \hat{s}_t \leq 0 \quad (17)$$

### C. Hourly Scenario-Check Subproblem

The hourly scenario check sub-problem (18) investigates the possibility of the operating point, $P_{it}$ and $I_{it}$, suggested by deterministic SCUC to be suitable in all probable scenarios [8]. Preventing actions modify the previous operating point to meet the constraints. The decision variables of this sub-problem are hourly power generation of thermal units, $P_{it}^c$. As mentioned in previous section, this sub-problem should be solved $2 \times m$ times for $m$ uncertain input variables. At the end, if the mean value of the objective values $S_t^c$, obtained from different runs of TPE, is larger than a predefined threshold, a Bender's Cut is referred to the master problem using (19). These violations are mostly due to insufficient dispatch capability of generating units in compensating deviations of probabilistic variables from their MPVs. Preventive actions removes all probable violations in some iterations.

$$\min S_t^c = s_t^c + s_{1t}^c + s_{2t}^c$$
$$s.t.$$

$$SF \cdot (K_P \cdot P_t^c - K_D \cdot P_{Dt}) - s_t^c \leq PL^{max}$$

$$-SF \cdot (K_P \cdot P_t^c - K_D \cdot P_{Dt}) - s_t^c \leq PL^{max}$$

$$\sum_{i=1}^{NG} P_{it}^c + \sum_{i=1}^{NW} P_{w,it}^c + s_{1t}^c - s_{2t}^c = \sum_{b=1}^{NB} P_{D,bt}^c$$

$$P_{it}^c - \hat{P}_{it} \leq R_i^{up} \cdot \hat{I}_{it} \quad ; \quad \lambda_{1,it}$$

$$\hat{P}_{it} - P_{it}^c \leq R_i^{dn} \cdot \hat{I}_{it} \quad ; \quad \lambda_{2,it}$$

$$P_{it}^c \leq P_i^{max} \cdot \hat{I}_{it} \quad ; \quad \mu_{1,it}$$

$$-P_{it}^c \leq -P_i^{min} \cdot \hat{I}_{it} \quad ; \quad \mu_{2,it}$$

$$s_t^c, s_{1t}^c, s_{2t}^c \geq 0 \qquad (18)$$

$$\sum_{i=1}^{NG} [(\hat{\lambda}_{1,it} - \hat{\lambda}_{2,it}) \cdot (P_{it} - \hat{P}_{it}) + (\hat{\lambda}_{1,it} \cdot R_i^{up} + \hat{\lambda}_{2,it} \cdot R_i^{dn} + \hat{\mu}_{1,it} \cdot P_i^{max} - \hat{\mu}_{2,it} \cdot P_i^{min}) \cdot (I_{it} - \hat{I}_{it})] + \hat{S}_t^c \leq 0 \qquad (19)$$

## V. CASE STUDIES

Two distinct cases are studied to illustrate the efficacy of the proposed algorithm. In the first case, a market based on six-bus system is investigated and the results of proposed algorithm is compared with the results of MCS as the benchmark and base case as the common approach. In the second case, the market is based on IEEE 118-bus system with a large number of uncertain variables which is not reported for this problem, so far: three wind farms and 91 probabilistic load.

### A. Six-Bus System

In the first case the proposed method is implemented on the six-bus system shown in Fig. 3. The system has three thermal generators, one wind farm, four transmission lines, two tap changers, and three load points. The parameters of the generation units, transmission system, hourly load, and hourly wind generation are given in [10]. Maximum dispatch adjustment capability of all thermal units are 10/60 of their ramp rate, because all corrective actions should be performed in ten minutes.

Three cases are studied in case A.

*A-1)* Deterministic SCUC with MPVs (as the base case)

*A-2)* Stochastic SCUC problem using MCS (as the benchmark)

*A-3)* Stochastic SCUC problem using TPE (the proposed method)

*Case A-1)* In the base case, uncertainties are modeled with MPVs. The base case scheduling is shown in Table I, and its total generation cost is 85 116 $.

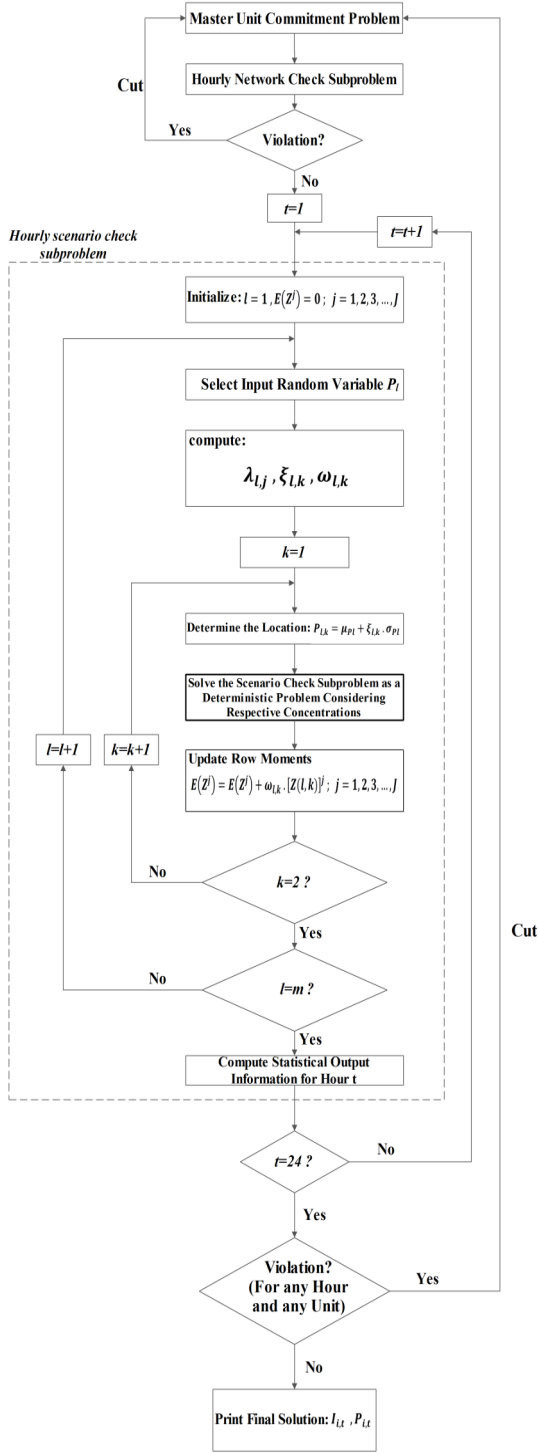

Fig. 2. Proposed method for Stochastic SCUC problem using TPE method

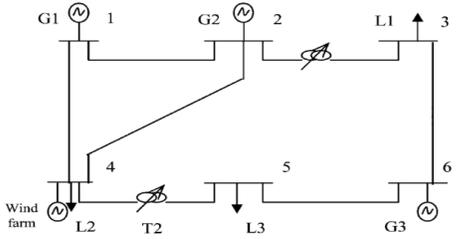

Fig. 3. One-line diagram of the six-bus system

TABLE I. Base case scheduling

| Unit | 1 | 2 | 3 |
|---|---|---|---|
| Commitment hours | 1-24 | 11-22 | 12-21 |

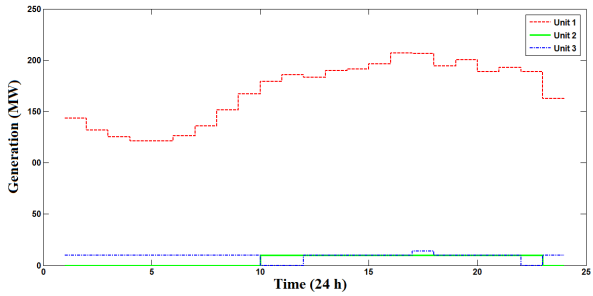

Fig. 4. Unit commitment and dispatch obtained by TPE method

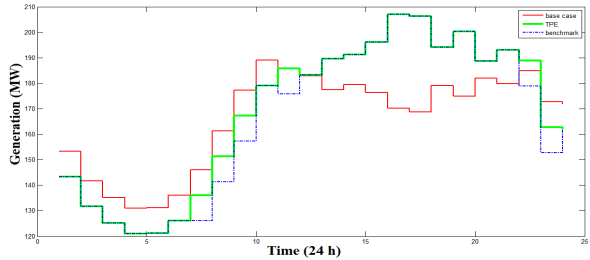

Fig. 5. Comparing generation schedules of unit 1 in three methods

TABLE II. Results of different approaches

|  | Case A-1 (Base Case) | Case A-2 (MCS) with SCENRED [8] Number of scenarios | | | Case A-3 (TPE) |
|---|---|---|---|---|---|
|  |  | S=10 | S=100 | S=1000 |  |
| Total generation cost ($) | 85166 | 87901 | 89891 | 90112 | 90014 |
| Relative Error (%) | 5.49% | 2.45% | 0.25% | 0% | 0.11% |
| Average CPU time (Sec) | 8 | 75 | 449 | 3114 | 70 |
| Number of iterations | 1 | 3 | 4 | 4 | 3 |

The MCS has been run three times with 10,100 and 1000 number of final scenarios, and the results are presented in Table II. Surely, the more number of final scenarios is, the more accuracy of the results is obtained. Therefore, the last case with 1000 final scenarios is considered as the benchmark to evaluate the proposed method. As the Table shows, for the benchmark, the algorithm converges in four iterations after 3114 seconds, and the total generation cost is 90112 $. It should be noticed that because the computer is not professional, the absolute CPU time is not judicable, but it is useful for comparing with the time of other approaches. The relative error is the relative difference between the total generation cost of each case and the benchmark.

*Case A-3)* In this case the Stochastic SCUC problem has been solved by the proposed method. Because there are four uncertain variables in the system, scenario-check sub-problems should be run eight times, in each iteration of the algorithm. The overall results and the schedule of thermal units are shown in Table II and Fig. 4, respectively. Fig. 5 compares the schedule of unit 1 in three different cases including the base case, the benchmark and the proposed approach. According to Table II, the proposed algorithm enjoys the advantage of low CPU time with acceptable accuracy. While the accuracy is much better than the MCS with 100 final scenarios, the CPU time is even less than the case of MCS with S=10. This is because of use of the concentrations instead of random scenarios. As Fig. 5 shows, the schedule of unit 1 in case 1-3 is very similar to the benchmark, so that in 17 hours of the day their schedules are the same while they are not so different in the other hours.

*Case A-2)* In this case the Stochastic SCUC problem has been solved by MCS as proposed in [8]. The hourly wind power is modeled with normal PDFs, and the nodal loads uncertainties are modeled with truncated normal PDFs with standard deviations of 10% of the mean values. The advantages of modeling of load points with truncated normal PDFs are explained in [7]. The standard deviation of the wind power is considered 20% of its mean values for each hour. To run the MCS, initially 10000 random scenarios are generated and then the scenarios are reduced into a number of final scenarios using fast backward/forward technique of SCENRED library of General Algebraic Modeling System (GAMS) software [18].

*B. IEEE 118-bus system*

The advantage of any speeding-up technique shows itself when the order of the system and the number of uncertain variables increase. In the case of scenario-based stochastic SCUC in large-scale power systems, the computational burden intensifies drastically, and efficient reduction of CPU time is very vital in this problem. To evaluate the proposed method in a large scale problem, the proposed approach is implemented on IEEE 118-bus system of [8] with a high number of uncertain variables.

The system has 54 thermal units, 186 transmission lines, 9 tap-changers and 91 load points [16]. The system is modified as presented in [19] by adding three wind farms W1, W2 W3, located on buses 36, 77 and 69, respectively. The power generated by the wind farms are given in [19]. The peak load of the system is 3733 MW and 12.78% of total consumed energy is procured by wind Farms.

The Stochastic SCUC problem is solved with two different number of uncertainties. In case B-1, only the uncertainty of three wind farms are considered (*m*=3), but in case B-2, in addition to the wind farms, the power consumption of all 91 load points are probable with normal PDFs (*m*=3+91).

The problem is solved once with MPVs and with deterministic SCUC approach, as the base case; and once with the proposed algorithm with TPE method. The total generation cost of the base case and TPE method are shown in Table III for both cases.

To evaluate the efficacy of proposed approach, in case A, MCS is administered with 1000 scenarios as the benchmark. In case B, due to the large scale of the network and large number of uncertainties, it is not reasonable to run the MCS with large number of scenarios. Instead, to appraise the success of proposed approach in this problem, similar to case A, the MCS is utilized through GAMS and 10000 final scenarios are selected. Then, the next day is simulated with the scenarios for both schedules resulted from the base case and the proposed approach.

To measure the ability of operating point to meet all probable variations of power generation or consumption; and also to determine the additional costs due to spinning of extra generators, two indices are proposed and formulated in the appendix A. The incapability of scheduled operating point to be corrected in real time by corrective actions is measured by the Corrective Actions Incapability (CAI) index. Extra Spinning Cost (ESC), is another index which indicates the extra operation cost imposed by spinning of newly-turned-on units because of preventive actions. ESC is the most possible waste of money due to change in the schedule from base case to the proposed method. Table IV shows CAI and ESC indices for both cases and with both techniques.

According to Table IV, CAI index of the base case is much higher than CAI of proposed approach. For instance, consider in case B-2 with 94 uncertain variables the market operator schedules the units according to the results of the base case. In the next day and in real time, 1.67% is the percentile probability that after variations in the wind or the loads, the operating reserves and the dispatch adjust capability of generating units are insufficient to compensate the power mismatch, and the system may enter into the alert state. Market scheduling using the proposed approach reduces the CAI to 0.2% at the expense of 3% increase of ESC. Surely, the last case is much more secure and suitable to be chosen by market operator.

VI. CONCLUSION

This paper presents a new algorithm for scenario-based Stochastic SCUC which uses Two Point Estimation (TPE) method to generate efficient concentrations instead of random scenarios for the algorithm. The algorithm enjoys the advantages of scenario-based techniques in such a complex problem, while significantly decreases the computational burden, which is the main disadvantage of these techniques. A lot of number of different uncertainties can be considered for Stochastic SCUC solution using TPE method. The results for the six-bus system and the modified IEEE 118-bus system illustrate the proposed method effectiveness. Implementing on the 118-bus system with 94 uncertain input variables has been shown that the proposed method can solve the Stochastic SCUC problem for large scale systems precisely with a reasonable simulation time.

APPENDIX

A. *Formulations of CAI and ESC indices*

TABLE III. Result of TPE approach

|  | Base case | TPE | |
| --- | --- | --- | --- |
|  |  | Number of uncertain input parameters | |
|  |  | m=3 | m=3+91 |
| Total cost ($) | 736011 | 743985 | 758369 |
| Average CPU time (min) | 30 | 100 | 240 |

Table IV. Evaluating the Corrective Actions Inability of the results

|  | Number of uncertain input parameters | Base case | TPE |
| --- | --- | --- | --- |
| CAI % | m=3 | 0.25% | 0.01% |
|  | m=3+91 | 1.67% | 0.20% |
| ESC % | m=3 | --- | 1.38% |
|  | m=3+91 | --- | 3.00% |

$$CAI = 1 - \sum_{S}^{NS}\sum_{i}^{NG} P_i^s \cdot A_i^s$$

$$A_i^s = \begin{cases} 1 & if\ (**)\ satisfied \\ 0 & others \end{cases}$$

$$R_i^{dn} \le \sum_{i=1}^{NG} P_{it}^s \cdot I_{it} + \sum_{i=1}^{NW} P_{w,it}^s - \sum_{b=1}^{NB} P_{D,bt}^s \le R_i^{up} \quad (**) \tag{A.1}$$

$$EC = \frac{\sum_{i \in x}\sum_{t \in x} C_i(t)}{Total\ Generation\ Cost}$$

$$x = \{(i,t)\ |\ I_{it}^{TPE} - I_{it}^{base\ case} = 1\};$$

$C_i(t)$: Constant Term of Production Cost Function  (A.2)

*B. Nomenclature*

**Variables**

| | |
|---|---|
| b | Index of buses. |
| $F_{ci}(.)$ | Production cost function of unit $i$. |
| i | Index of units. |
| $I_{it}$ | Commitment of unit $i$ at time $t$. |
| $P_{it}$ | Power dispatch of unit $i$ at time $t$. |
| $P_{D,t}$ | System demand at time $t$. |
| $P_{D,bt}$ | System demand at time $t$ at bus $b$. |
| $P_{w,it}$ | Generation of wind power unit $i$ at time $t$. |
| $R_{o,it}$ | Operating reserve of unit $i$ at time $t$ |
| $R_{s,it}$ | Spinning reserve of unit $i$ at time $t$ |
| $S^{(.)}_{(.)}$ | Slack variable. |
| $SU_{it}$, $SD_{it}$ | Startup/shutdown cost of unit $i$ at time $t$. |
| t | Index of hours. |
| $X_{it}^{off}$ | Off time of unit $i$ at time $t$ |
| $X_{it}^{on}$ | On time of unit $i$ at time $t$ |
| $\lambda$, $\mu$ | Dual variables. |
| $(.)^C$ | Variable related to concentration $C$. |
| $(.)^f$ | Indicating forecasted value. |
| $(.)^S$ | Variable related to scenario $S$. |

**Constants**

| | |
|---|---|
| $DR_i$ | Ramp-down rate limit of unit $i$ |
| NB | Number of buses. |
| NG | Number of non-wind units. |
| NT | Number of periods under study (24 h). |
| NW | Number of wind power units. |
| $P_i^{max}$ | Maximum capacity of unit $i$. |
| $P_i^{min}$ | Minimum capacity of unit $i$. |
| $R_i^{up}$, $R_i^{dn}$ | UP/down limits for corrective dispatch |
| $R_{o,t}$ | System operating reserve requirement |
| $R_{s,t}$ | System spinning reserve requirement |
| $T_i^{off}$ | Minimum off time of unit $i$ |
| $T_i^{on}$ | Minimum on time of unit $i$ |
| $UR_i$ | Ramp-up rate limit of unit $i$ |

**Matrices**

| | |
|---|---|
| $K_P$, $K_D$ | Bus-generator/bus-load incidence matrix. |
| $P_{Dt}$ | Vector of system demand at time $t$. |
| $P_t$ | Vector of generating power at time $t$. |
| $PL^{max}$ | Vector of upper limit for power flow. |